\pdfoutput=1
\documentclass{amsart}
\usepackage{latexsym}
\usepackage{amsfonts}
\usepackage{graphicx}
\newcommand{\C}{\mathbb{C}}
\newcommand{\R}{\mathbb{R}}
\newcommand{\To}{\rightarrow}

\begin{document}
\title{The Asymptotic Variety of a Pinchuk Map as a Polynomial Curve}
 \author{L. Andrew Campbell}
\address{908 Fire Dance Lane \\
Palm Desert CA 92211 \\ USA}
\email{landrewcampbell@earthlink.net}

\begin{abstract}
The asymptotic variety of a counterexample
of Pinchuk type  to the
strong real Jacobian conjecture is
explicitly described by low degree polynomials.
\end{abstract}

\maketitle
\section{Introduction}
Let the polynomial map $F=(P,Q): \R^2 \To \R^2$
be the Pinchuk map of total degree $25$ considered in
\cite{PPR,ArnoBook,Picturing,PPRErr}.
Its asymptotic variety, $A(F)$, a closed curve in the
image $(P,Q)$-plane,
was computed in \cite{PPR}. It is depicted below using
differently scaled axes. It intersects the vertical axis at $(0,0)$
and $(0,208)$ and its leftmost point is $(-1,-163/4)$.

\begin{figure}[ht]
\centerline{
\includegraphics[height=2in,width=6in]{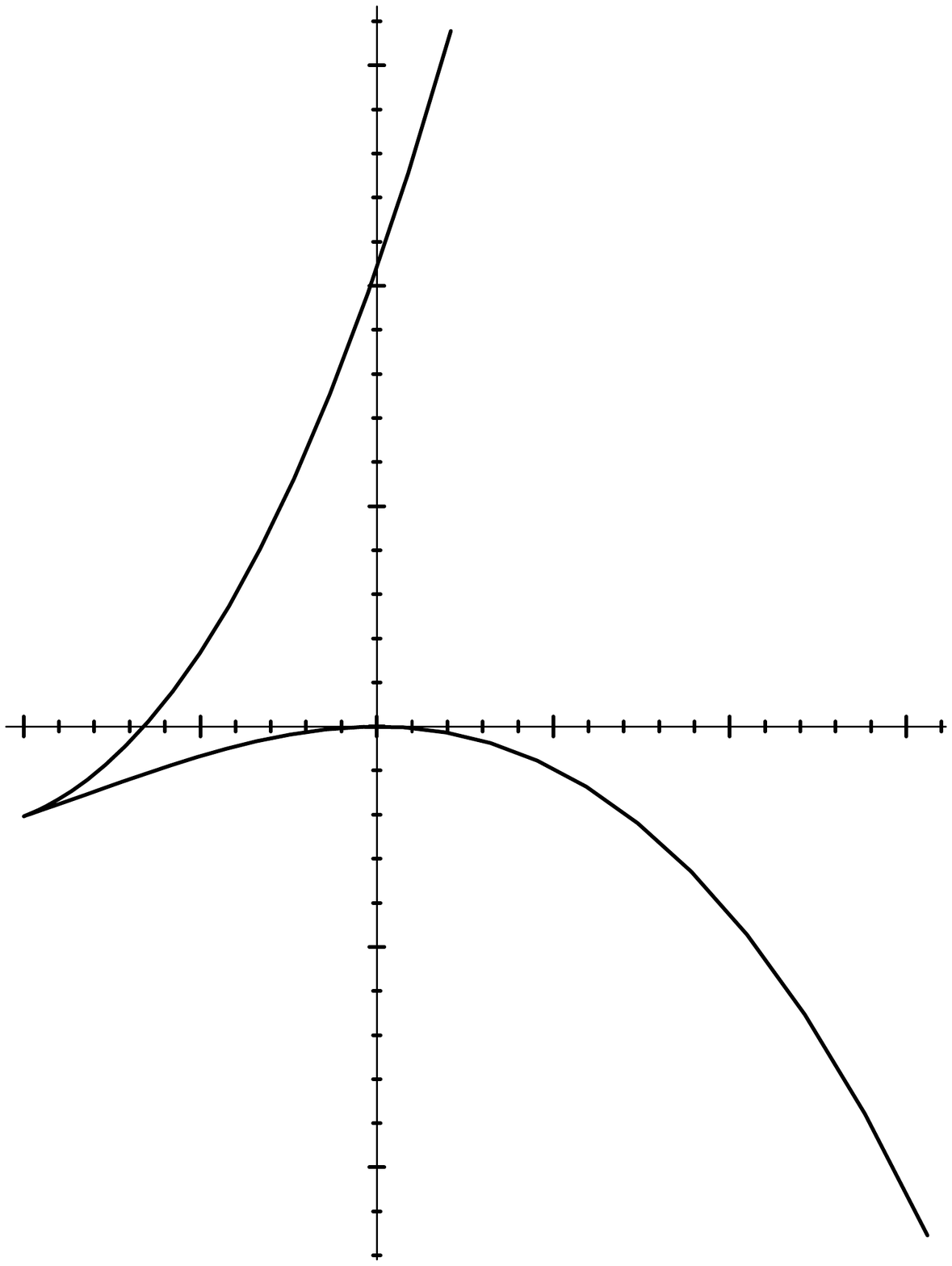}}
\caption{The asymptotic variety of the Pinchuk map $F$.}
\end{figure}

This brief note will show that $A(F)$ has the bijective polynomial parametrization by  $s \in \R$:
\begin{equation}
P(s) = s^2 - 1
\end{equation}
\begin{equation}
Q(s) = -75s^5 +
\frac{345}{4}s^4 - 29s^3 + \frac{117}{2}s^2
 - \frac{163}{4}
\end{equation}
and that its points satisfy the minimal equation
\begin{equation}
(Q -(345/4)P^2 -231P -104)^2 = (P+1)^3(75P + 104)^2
\end{equation}
In particular, only one point on the curve satisfies $P = -1$,
and that point is its only singular point. Also, there is a single point with $P=-104/75$ which satisfies Eqn (3). That point
is not on the curve itself, because Eqn (1) implies $P \ge -1$,
but belongs to its Zariski closure.

Exactly parallel facts were established in \cite{GeoPinMap}
for a different Pinchuk map $\tilde F$ of total degree $40$. In fact, any two Pinchuk maps have essentially the same behavior and asymptotic variety, differing only by a triangular polynomial automorphism of the image plane.

A parametrization of $A(F)$ appeared in the unpublished preprint \cite{Picturing}. The derivation is shortened here. $F$ is a useful  reference example, because of the simplicity and low degree of the explicit equations for $A(F)$.

\section{Pinchuk maps}
Pinchuk maps are certain polynomial maps $F=(P,Q): \R^2 \to \R^2$
that have an everywhere positive Jacobian determinant $j(P,Q)$,
and are not injective \cite{Pinchuk}. The polynomial $P(x,y)$
is constructed by defining
$t=xy-1,h=t(xt+1),f=(xt+1)^2(t^2+y),P=f+h$.
Note that $\deg h = 5$, $\deg f = 10$, so $\deg P = 10$.
The polynomial $Q$ varies for different Pinchuk maps,
but always has the form
$Q = -t^2 -6th(h+1) - u(f,h)$,
where $u$ is an auxiliary polynomial in $f$ and $h$,
chosen so that
$j(P,Q) = t^2 + (t+f(13+15h))^2 + f^2$.
As in \cite{PPR,ArnoBook}, choose specifically
\begin{equation}
u =
170fh + 91h^2 + 195fh^2 + 69h^3 + 75fh^3+
\frac{75}{4}h^4.
\end{equation}
Then $\deg F = \deg Q = 25$.

Suppose $\tilde F = (P,\tilde Q)$ is a different Pinchuk map
defined using $\tilde u$. Observe that $A = Q - \tilde Q = u(f,h) -
\tilde u(f,h)$ lies in $\R[P,h]$ and satisfies $j(P,A) = 0$, since $j(P,\tilde Q) = j(P,Q)$. By \cite{DedekindD} or \cite[Thm. 1.2.25]{ArnoBook}, the subalgebra $\R[P,A] \subset \R[x,y]$ is generated by a single element. In $\R[P,h]$ that  generator must be
a degree $1$ polynomial in $P$ alone, since $P$ and $h$ are algebraically independent. So $A \in R[P]$ and  $\tilde F = T \circ F$ for a triangular polynomial automorphism $T(x,y) = (x,y+S(x))$.

The ``original'' Pinchuk map of \cite{GeoPinMap} is defined by adding, not subtracting, $ (1/4)f(75f^3+300f^2h+450fh^2 +276f^2
          +828fh+48h^2+364f +48h)$, which is thus $- \tilde u$. Clearly $\deg \tilde F = \deg \tilde Q = 40$.

Note that no Pinchuk map can have degree less than $25$.
For $P$ has degree $10$, and so, if
$\tilde Q = Q + S(P)$, there is no way to cancel the terms
of degree $25$ without introducing terms of yet higher degree.

\section{Asymptotic behavior}
The points $(-1,-163/4)$ and $(0,0)$ of $A(F)$ have no inverse image under $F$, all other points of $A(F)$ have one inverse
image, and all points of the image plane not on $A(F)$ have two.
See \cite{PPR,PPRErr}.

If the two omitted points are deleted from $A(F)$, we are left
with three curves. Since the curves tend to infinity or to an omitted point at either end, their inverse images tend to infinity at both ends. So they partition their complement, $\R^2 \setminus F^{-1}(A(F))$,
into four simply connected domains.
These domains are mapped homeomorphically to their images, two each to the domains on either side of $A(F)$.
See \cite{Picturing}.

Suppose $\tilde F = (P,\tilde Q) = T \circ F$ is a different Pinchuk map, with asymptotic variety $A(\tilde F)$. From the definition of the asymptotic variety of a polynomial map as the set of finite limits of the map along curves that tend to infinity (\cite{AsympVals,Asymptotics}), or equivalently,  the points at which the map is not proper (\cite{notproper,ZJGeoRealMaps}), it follows that $A(\tilde F) = T(A(F))$. The behavior is that of $F$, up to
the triangular automorphism $T$ of the image plane.

Also, note that  the partition of the $(x,y)$-plane  into three curves and four domains is exactly the same as for $F$.. See \cite{GeoPinMap} for a graphic depiction.

\section{Equations for $A(F)$}
Also from previously cited work, a general level set $P=c$
in the $(x,y)$-plane has a rational parametrization 
$$x(h) = \frac{ (c-h)(h+1) }{  (c-2h-h^2)^2 }$$
$$y(h) = \frac{ (c-2h-h^2)^2(c-h-h^2) }{ (c-h)^2 },$$
which can be obtained by solving $P = c$ for $x$ and then
$y$, and can be readily verified by substitution into the
defining equations $t=xy-1,h=t(xt+1),f=(xt+1)^2(t^2+y),P=f+h$.
Note that $h(x(h),y(h))$ does indeed simplify to just $h$,
and $P(x(h),y(h))$ to just $c$.

Temporarily ignore the special cases $c=-1$ and $c=0$, for which different parametrizations apply. At poles
$c=h$, it is easy to check that $Q = -t^2 -6th(h+1) - u(f,h)$
is infinite. Indeed, $f=c-h=0$ and the dominant pole is
$-h^4(h+1)^2/(c-h)^2$ with a nonzero numerator by the
restriction on $c$.
In contrast, as $h$ tends to a pole $c=h^2+2h$ we find that $Q$ approaches a finite limit. The limit  condition
on $c$ and $h$ can also be stated as $h = -1 \pm \sqrt{1+c}$.
 This yields two points of $A(F)$ on a vertical line $P=c$
 when $c > -1$ and $ c \ne 0$, and no points when $c < -1$.

In more detail, $x(h)y(h) =(h+1)(c-h-h^2)(c-h)^{-1}$.
If we substitute $c=h^2+2h$, the expression simplifies to
$h(h+1)h^{-1}(h+1)^{-1}$. Here too, $h \ne 0$ and $h \ne -1$
by the restriction on $c$, so the ratio is defined and equal to $1$.
This only means that as $h$ tends to such a pole, $xy$ tends
to $1$.
In the limit $t=0, f =c - h = h^2+h$ and thus
$Q =  -t^2 -6th(h+1) - u(f,h) =-u(h^2+h,h)$.
Since $A(F)$ is closed in the Euclidean topology, it contains the three sofar missing points with $P=-1$ or  $P=0$ and the same equations hold there by continuity.
That is all of $A(F)$, because every point not obtained is, by section 3, known not to be an  asymptotic value of $F$.

The polynomial parametrization
$(P(h),Q(h)) = (h^2+2h,-u(h^2+h,h)$
is clearly a bijection from $\R$ onto $A(F)$.
This works (with the appropriate $u$)
for any Pinchuk map; it is the form
reported in \cite{GeoPinMap}.

Using $s = h +1$ as a parameter instead and the specific
auxiliary polynomial $u(f,h)$ in Eqn (4) yields Eqn (1)
and Eqn (2). The choice of $s$ simplifies the calculation
of the gradient of the parametrization and of points on $A(F)$.
From Eqn (2)
$$
Q - \frac{345}{4}s^4 - \frac{117}{2}s^2 + \frac{163}{4}
=
- s(75s^4 +  29s^2).
$$
Squaring both sides and substituting $P+1$ for $s^2$
yields Eqn (3). Expand and rewrite Eqn (3) as an implicit polynomial equation  $B(P,Q)=0$, quadratic in $Q$. $B$ cannot have a factor
that is a nonconstant polynomial in $P$, because the coefficient of $Q^2$ in $B$ is $1$. Nor can $B$ have two factors linear in $Q$. At least one such factor, say $Q - K(P)$ for a polynomial $K$, would have to be identically $0$ on $A(F)$, yet some  vertical lines $P=c$ do not intersect $A(F)$. So $B$ is irreducible and therefore its set of zeroes is the Zariski closure of $A(F)$.

\section{Double asymptotic identities}
Ronen Peretz championed a simpler way of finding parametrization
equations such as Eqn (1) and Eqn (2) \cite{AsympVals,Asymptotics}.
A double asymptotic identity for $F$ is an equation
$F(R(x,y)) = G(x,y)$ for a  rational (but not polynomial)
map $R$ and a polynomial map $G$.

Consider $R = (x^{-2},yx^3+x^2)$. As $t \circ R = xy,
h \circ R = (x+y)y,  f \circ R = (x+y)^2(y^2+xy+1)$,
the map  $G = (P \circ R,Q \circ R)$ is polynomial.

As $x$ tends to  zero for a fixed $y$, the point $R(x,y)$ tends
to infinity, describing a curve along which $F$ tends to the finite limit point $G(0,y)$. So $G(0,y) = (y^4+2y^2,-u(y^4+y^2,y^2))$
is a parametrization of (some of the points of) $A(F)$.

Comparing with the bijective parametrization by $h$ of the previous section, it is evident that this parametrization  covers only the points $h \ge 0$. Each such point is obtained twice, except for $(P,Q)=(0,0)$, where the parametrization reverses course. A similar parametrization covering the points $h \le 0$ of $A(F)$ arises from the alternate choice of $(-x^{-2},yx^3-x^2)$ for $R$ \cite{Picturing}.

The computations of the previous section can be recast into
a rational  identity of the form $ F(R(x,y)) = G(x,y)$ that
provides the  bijective parametrization. However, $G$ is not a
polynomial map, but rather a rational map with
$G(0,y)$ defined for all but finitely many values of $y$
and polynomial in $y$. The two exceptional values of $y$
correspond to the special cases $P=-1$ and $P=0$.

\section{Relation to the Jacobian conjecture}
A weak Jacobian conjecture for polynomial maps of $\R^2$
to itself is that a Keller map (nonzero constant Jacobian determinant)
is injective. This is weaker than the standard Jacobian conjecture JC(2,$\R$), even though injectivity implies bijectivity
 here, because the inverse is not required to be a polynomial map.

{\em Remark}. Over $\C$ this distinction does not exist, since
any inverse map is birational and everywhere defined, hence
polynomial. Both conjectures for $\R^2$ would follow from
JC(2,$\C$). Note that JC(2,$\R$) is not known to imply
JC(2,$\C$).

A general feature of the Pinchuk map $F$ that conflicts with the Keller condition is that radial similarity of Newton polygons fails.
$N(P) = N(x^6y^4 + x^2 + y)$, a quadrilateral, while
$N(Q) = N(x^{15}y^{10} + x^3y^4 + x^5 + y)$, a five-sided
polygon. For $\tilde F$, though,
$N(\tilde Q) = N(x^{24}y^{16} + x^8 + y^4)$, a 
$4$-fold radial expansion of $N(P)$. All these polygons have no edge of negative slope.

For a nonsingular map $f = (p,q):\R^2 \To \R^2$, whether
polynomial or not, $j(p,q)$ is the rate of change of $q$ along the level curves of $p$, parametrized as the flow of the Hamiltonian vector field
$H(p) = (-\partial p/\partial y , \partial p/\partial x)$. If $f$ is polynomial the local flow $(x(t),y(t))$ with initial condition $x(0) = a,  y(0) = b$ is at least real analytic. If $f$ is also a Keller map, it has a polynomial inverse if the power series expansions for
$x$ and $y$ have infinite radius of convergence for even a single point $(a,b) \in \R^2$, in which case the flow is actually polynomial for any $(a,b) \in \R^2$. That follows from the corresponding result for complex Keller maps \cite[Thm 3.2]{HamiltonianFlows}, by treating $f$ as a complex polynomial map and $(a,b)$ as a point of $\C^2$.

\section{Acknowledgments}
Eqn (4) was first circulated by Arno van den Essen in an
email message to a number of colleagues in June 1994.
I made two blunders and a significant typographical error in
describing the associated Pinchuk  map $F$ in \cite{PPR}. The errors were  corrected
 in \cite{Picturing} and belatedly for the official record
in \cite{PPRErr}. Janusz Gwo{\'z}dziewicz significantly clarified
the relationship between different Pinchuk maps for me.

\end{document}